\documentclass[preprint,12pt]{elsarticle}

\usepackage{amssymb}

\usepackage{algorithm}
\usepackage{algorithmic}

\newtheorem{definition}{Definition}

\newtheorem{lemma}{Lemma}
\newtheorem{remark}{Remark}
\newtheorem{theorem}{Theorem}

\journal{Pattern Recognition}

\begin{document}

\begin{frontmatter}



\title{Iterative thresholding algorithm based on non-convex method for modified $\l_{p}$-norm regularization minimization}

 \author[label1]{Angang Cui}
 \ead{cuiangang@163.com}
 \author[label2]{Jigen Peng \corref{cor}}
 \cortext[cor]{Corresponding author}
 \ead{jgpengxjtu@126.com}
  \author[label3]{Haiyang Li}
  \ead{fplihaiyang@126.com}
  \author[label3]{Meng Wen}
 \ead{wen5495688@163.com}
 \author[label1]{Junxiong Jia}
 \ead{jjx323@xjtu.edu.cn}
 \address[label1]{School of Mathematics and Statistics, Xi'an Jiaotong University, Xi'an, 710049, China.}
   \address[label2]{School of Mathematics and Information Science, Guangzhou University, Guangzhou, 510006, China.}
 \address[label3]{School of Science, Xi'an Polytechnic University, Xi'an, 710048, China.}

\begin{abstract}
Recently, the $\l_{p}$-norm regularization minimization problem $(P_{p}^{\lambda})$ has attracted great attention in compressed sensing. However, the $\l_{p}$-norm $\|x\|_{p}^{p}$ in
problem $(P_{p}^{\lambda})$ is nonconvex and non-Lipschitz for all $p\in(0,1)$, and there are not many optimization theories and methods are proposed to solve this problem. In fact,
it is NP-hard for all $p\in(0,1)$ and $\lambda>0$. In this paper, we study two modified $\l_{p}$ regularization minimization problems to approximate the NP-hard problem $(P_{p}^{\lambda})$.
Inspired by the good performance of Half algorithm and $2/3$ algorithm in some sparse signal recovery problems, two iterative thresholding algorithms are proposed to solve the problems $(P_{p,1/2,\epsilon}^{\lambda})$ and $(P_{p,2/3,\epsilon}^{\lambda})$ respectively. Numerical results show that our algorithms perform effectively in finding the sparse signal in some
sparse signal recovery problems for some proper $p\in(0,1)$.
\end{abstract}

\begin{keyword}
Compressed sensing\sep Modified $\l_{p}$ regularization minimization problem\sep Thresholding representation theory\sep $1/2-\epsilon$ algorithm\sep $2/3-\epsilon$ algorithm

\MSC 90C26\sep 65K10\sep 49M20

\end{keyword}

\end{frontmatter}


\section{Introduction}\label{section1}
During the last decade, the $\l_{p}$-norm regularization minimization problem has attracted great attention in compressed sensing \cite{chen1,xu2,cao3,chen4,chen5}. In mathematics,
it can be modeled into the following minimization problem
\begin{equation}\label{equ1}
(P_{p}^{\lambda})\ \ \ \ \ \min_{z\in \mathbb{R}^{n}}\Big\{\|Az-b\|_{2}^{2}+\lambda \|z\|_{p}^{p}\Big\}
\end{equation}
for some $\lambda\in(0,+\infty)$ and $p\in(0,1)$, where $A\in \mathbb{R}^{m\times n}$ is a real matrix of full row rank with $m\ll n$, $b\in \mathbb{R}^{m}$ is a nonzero real column
vector and $\|z\|_{p}^{p}=\sum_{i=1}^{n}|z_{i}|^{p}$
for any $z\in \mathbb{R}^{n}$. We can see that as $p\rightarrow 0$, the problem $(P_{p}^{\lambda})$ tends to the $\l_{0}$-norm regularization minimization problem
\begin{equation}\label{equ2}
(P_{0}^{\lambda})\ \ \ \ \ \min_{z\in \mathbb{R}^{n}}\Big\{\|Az-b\|_{2}^{2}+\lambda \|z\|_{0}\Big\},
\end{equation}
where $\|z\|_{0}$ denotes the number of nonzero components of $z$, and the iterative hard thresholding algorithm (Hard algorithm)\cite{thom6} has been proposed for solving the
$\l_{0}$ regularization minimization problem $(P_{0}^{\lambda})$. In addition, as $p\rightarrow 1$, the problem $(P_{p}^{\lambda})$ tends to the $\l_{1}$-norm regularization
minimization problem
\begin{equation}\label{equ3}
(P_{1}^{\lambda})\ \ \ \ \ \min_{z\in \mathbb{R}^{n}}\Big\{\|Az-b\|_{2}^{2}+\lambda \|z\|_{1}\Big\},
\end{equation}
where $\|z\|_{1}=\sum_{i=1}^{n}|z_{i}|$, and $z_{i}$ represents the $i$-th component of vector $z$. As the compact convex relaxation of the NP-hard problem $(P_{0}^{\lambda})$,
there are many efficient methods (e.g., see \cite{dau7,don8,tom9,yin10,yang11}) being proposed for solving the $\l_{1}$-norm regularization minimization problem $(P_{1}^{\lambda})$.
It is obvious that the problem $(P_{p}^{\lambda})$ is intermediate between the $\l_{0}$-norm regularization minimization problem $(P_{0}^{\lambda})$ and the $\l_{1}$-norm
regularization minimization problem $(P_{1}^{\lambda})$ for any $p\in(0,1)$, because of the relationship
\begin{equation}\label{equ4}
\lim_{p\rightarrow 0^{+}}\sum_{i=1, z_{i}\neq 0}^{n}|z_{i}|^{p}=\|z\|_{0}\ \ \ \mathrm{and} \ \ \ \lim_{p\rightarrow 1}\sum_{i=1 0}^{n}|z_{i}|^{p}=\|z\|_{1}.
\end{equation}

Unfortunately, the  problem $(P_{p}^{\lambda})$ is a nonconvex and non-Lipschitz minimization problem. There are not many optimization theories on analyzing this type of
problems and it is NP-hard for all $p\in(0,1)$ and $\lambda>0$ (see \cite{chen1}). At present, the most direct way to solve the problem $(P_{p}^{\lambda})$ is that the
iterative thresholding algorithm only when $p=1/2, 2/3$ (see \cite{xu2,cao3}). In fact, the corresponding thresholding functions for the problem $(P_{p}^{\lambda})$
are in closed form only for $p=1/2, 2/3$. Xu et al.\cite{xu2} and Cao et al.\cite{cao3} have shown that the problem $(P_{p}^{\lambda})$ could be fast solved by the iterative
$\l_{1/2}$ thresholding algorithm (Half algorithm) and iterative $\l_{2/3}$ thresholding algorithm ($2/3$ algorithm), and the computational complexity of these two iterative
algorithms are all $\mathcal{O}(mn)$. A major drawback of the iterative thresholding algorithm for the problem $(P_{p}^{\lambda})$ is that the closed form iterative
thresholding for the problem $(P_{p}^{\lambda})$ available only at $p=1/2,2/3$.

In this paper, we propose a modified $\ell_{p}$-norm to replace the nonconvex and non-Lipschitz  $\l_{p}$ norm $\|x\|_{p}^{p}$ given by
\begin{equation}\label{equ5}
\sum_{i=1}^{n}\frac{|z_{i}|^{\theta}}{(|z_{i}|+\epsilon_{i})^{\theta-p}}
\end{equation}
for all $p\in(0,1)$, where $\theta>0$ and $\epsilon_{i}>0$. With the change of parameter $\epsilon_{i}>0$, we have
\begin{equation}\label{equ6}
\lim_{\epsilon_{i}\rightarrow0^{+}}\frac{|z_{i}|^{\theta}}{(|z_{i}|+\epsilon_{i})^{\theta-p}}=|z_{i}|^{p},
\end{equation}
and the modified $\ell_{p}$-norm (\ref{equ6}) approximates the $\l_{p}$-norm of vector $z$:
\begin{equation}\label{equ7}
\lim_{\epsilon_{i}\rightarrow0^{+}}\sum_{i=1}^{n}\frac{|z_{i}|^{\theta}}{(|z_{i}|+\epsilon_{i})^{\theta-p}}=\|z\|^{p}_{p}.
\end{equation}
Therefore, the $\l_{p}$-norm regularization minimization problem transformed by the modified $\ell_{p}$-norm (\ref{equ6}) could be written as the following minimization problem
\begin{equation}\label{equ8}
(P_{p,\theta,\epsilon}^{\lambda})\ \ \ \ \ \min_{z\in \mathbb{R}^{n}}\Big\{\|Az-b\|_{2}^{2}+\lambda \sum_{i=1}^{n}\frac{|z_{i}|^{\theta}}{(|z_{i}|+\epsilon_{i})^{\theta-p}}\Big\}
\end{equation}
for all $p\in(0,1)$. In particular, we do claim that the problem $(P_{p,\theta,\epsilon}^{\lambda})$ matches the following special version
\begin{equation}\label{equ9}
(P_{p,1/2,\epsilon}^{\lambda})\ \ \ \ \ \min_{z\in \mathbb{R}^{n}}\Big\{\|Az-b\|_{2}^{2}+\lambda \sum_{i=1}^{n}\frac{|z_{i}|^{1/2}}{(|z_{i}|+\epsilon_{i})^{1/2-p}}\Big\}
\end{equation}
for $\theta=1/2$, and
\begin{equation}\label{equ10}
(P_{p,2/3,\epsilon}^{\lambda})\ \ \ \ \ \min_{z\in \mathbb{R}^{n}}\Big\{\|Az-b\|_{2}^{2}+\lambda \sum_{i=1}^{n}\frac{|z_{i}|^{2/3}}{(|z_{i}|+\epsilon_{i})^{2/3-p}}\Big\}
\end{equation}
for $\theta=2/3$.

Throughout this paper, we just consider the special problems $(P_{p,1/2,\epsilon}^{\lambda})$ and $(P_{p,2/3,\epsilon}^{\lambda})$ for all $p\in(0,1)$, and extend the aforementioned
well-known Half algorithm and the $2/3$ algorithm \cite{xu2,cao3} to solve these two problems. The outline of this paper is as follows. In Section \ref{section2}, some preliminary
results used in this paper are given. In Section \ref{section3}, we propose two iterative thrsholding algorithms to solve the problems $(P_{p,1/2,\epsilon}^{\lambda})$ and
$(P_{p,2/3,\epsilon}^{\lambda})$ respectively. In Section \ref{section4}, we conduct some numerical experiments to show the performance of our algorithm. Some conclusion remarks are
presented in Section \ref{section5}.

\section{Preliminaries} \label{section2}
In this section, we give some crucial preliminary results that are used in this paper.

\begin{lemma}\label{lem1}{\rm(see \cite{xu2})}
For any fixed $\lambda>0$ and $\beta,r\in \mathbb{R}$, suppose that
\begin{equation}\label{equ11}
h_{1/2,\lambda}(r)\triangleq\arg\min_{\beta\in \mathbb{R}}\Big\{(\beta-r)^{2}+\lambda|\beta|^{1/2}\Big\},
\end{equation}
then the operator $h_{1/2,\lambda}(r)$ can be expressed by
\begin{equation}\label{equ12}
h_{1/2,\lambda}(r)=\left\{
    \begin{array}{ll}
      f_{1/2,\lambda}(r), & \ \mathrm{if} \ {|r|> \frac{\sqrt[3]{54}}{4}\lambda^{2/3};} \\
      0, & \ \mathrm{if} \ {|r|\leq\frac{\sqrt[3]{54}}{4}\lambda^{2/3}.}
    \end{array}
  \right.
\end{equation}
where
\begin{equation}\label{equ13}
f_{1/2,\lambda}(r)=\frac{2}{3}r\Big(1+\cos\Big(\frac{2\pi}{3}-\frac{2}{3}\arccos\Big(\frac{\lambda}{8}\Big(\frac{|r|}{3}\Big)^{-3/2}\Big)\Big)\Big).
\end{equation}
\end{lemma}

\begin{lemma}\label{lem2}{\rm(see \cite{cao3})}
For any fixed $\lambda>0$ and $\beta,r\in \mathbb{R}$, suppose that
\begin{equation}\label{equ14}
h_{2/3,\lambda}(r)\triangleq\arg\min_{\beta\in \mathbb{R}}\Big\{(\beta-r)^{2}+\lambda|\beta|^{2/3}\Big\},
\end{equation}
then the operator $h_{2/3,\lambda}(r)$ can be expressed by
\begin{equation}\label{equ15}
h_{2/3,\lambda}(r)=\left\{
    \begin{array}{ll}
      f_{2/3,\lambda}(r), & \ \mathrm{if} \ {|r|> \frac{\sqrt[4]{48}}{3}\lambda^{3/4};} \\
      0, & \ \mathrm{if} \ {|r|\leq\frac{\sqrt[4]{48}}{4}\lambda^{3/4}.}
    \end{array}
  \right.
\end{equation}
where
\begin{equation}\label{equ16}
f_{2/3,\lambda}(r)=\frac{1}{8}\Bigg(|\Phi_{2/3,\lambda}(r)|+\sqrt{\frac{2|r|}{|\Phi_{2/3,\lambda}(r)|}-|\Phi_{2/3,\lambda}(r)|^{2}}\Bigg)^{3}sign(r)
\end{equation}
and
\begin{equation}\label{equ17}
\Phi_{2/3,\lambda}(r)=\frac{2}{\sqrt{3}}\lambda^{1/4}\Big(\cosh\Big(\frac{1}{3}\mathrm{arccosh}\Big(\frac{27}{16}\lambda^{-3/2}r^{2}\Big)\Big)\Big)^{1/2}.
\end{equation}
\end{lemma}

\begin{definition}\label{def1}
(\cite{sim12}) The nonincreasing rearrangement of the vector $x\in \mathbb{R}^{n}$ is the vector $\lceil x\rfloor\in \mathbb{R}^{n}$ for which
$$\lceil x\rfloor_{1}\geq \lceil x\rfloor_{2}\geq\cdots\geq \lceil x\rfloor_{n}\geq0$$
and there is a permutation $\pi:[n]\rightarrow[n]$ with $\lceil x\rfloor_{i}=|x_{\pi(i)}|$ for all $i\in[n]$.
\end{definition}

\section{Two iterative thresholding algorithms for solving problem $(P_{p,\theta,\epsilon}^{\lambda})$}\label{section3}
In this section, we propose two iterative thresholding algorithms, namely, $1/2-\epsilon$ algorithm ($\theta=1/2$) and the $2/3-\epsilon$ algorithm ($\theta=2/3$), to
solve the problems $(P_{p,1/2,\epsilon}^{\lambda})$ and $(P_{p,2/3,\epsilon}^{\lambda})$ respectively. Moreover, we also provide some convergence analysis for our methods.
We should declare that the study of the $1/2-\epsilon$ algorithm and the $2/3-\epsilon$ algorithm proposed below are motivated by the well-known Half algorithm and $2/3$
algorithm proposed in Xu et al. \cite{xu2} and Cao et al. \cite{cao3}.

\subsection{The $1/2-\epsilon$ algorithm for solving the problem $(P_{p,1/2,\epsilon}^{\lambda})$}\label{subsection3-1}
In the subsection, we propose the $1/2-\epsilon$ algorithm to solve the problem $(P_{p,1/2,\epsilon}^{\lambda})$ for all $p\in(0,1)$.
Before the analytic expression of the $1/2-\epsilon$ algorithm, we should derive the closed form representation of the optimal solution to the problem
$(P_{p,1/2,\epsilon}^{\lambda})$, which underlies the algorithm to be proposed.

For any $\lambda,\mu\in(0,\infty)$, $p\in(0,1)$ and $z, y\in \mathbb{R}^{n}$, let
\begin{equation}\label{equ18}
\mathcal{C}^{1/2}_{\lambda}(z)=\|Az-b\|_{2}^{2}+\lambda \sum_{i=1}^{n}\frac{|z_{i}|^{1/2}}{(|z_{i}|+\epsilon_{i})^{1/2-p}},
\end{equation}
\begin{equation}\label{equ19}
\begin{array}{llll}
\mathcal{C}_{\lambda,\mu}^{1/2}(z,y)&=&\mu\|Az-b\|_{2}^{2}+\displaystyle\lambda\mu \sum_{i=1}^{n}\frac{|z_{i}|^{1/2}}{(|y_{i}|+\epsilon_{i})^{1/2-p}}-\mu\|Az-Ay\|_{2}^{2}+\|z-y\|_{2}^{2}
\end{array}
\end{equation}
and
\begin{equation}\label{equ20}
B_{\mu}(z)=z+\mu A^{\top}(b-Az).
\end{equation}

\begin{lemma}\label{lem3}
For any $\lambda,\mu\in(0,\infty)$ and $p\in(0,1)$, if $\tilde{z}=(\tilde{z}_{1},\tilde{z}_{2},\cdots, \tilde{z}_{n})^{\top}$ is  a local minimizer to $\mathcal{C}_{\lambda, \mu}^{1/2}(z,y)$, then
\begin{equation}\label{equ21}
\tilde{z}_{i}=0\Leftrightarrow |[B_{\mu}(y)]_{i}|\leq t_{1/2,\lambda\mu/(|y_{i}|+\epsilon_{i})^{1/2-p}}
\end{equation}
and
\begin{equation}\label{equ22}
\tilde{z}_{i}=f_{1/2,\lambda\mu/(|y_{i}|+\epsilon_{i})^{1/2-p}}([B_{\mu}(y)]_{i})\Leftrightarrow |[B_{\mu}(y)]_{i}|> t_{1/2,\lambda\mu/(|y_{i}|+\epsilon_{i})^{1/2-p}},
\end{equation}
where $[B_{\mu}(y)]_{i}$ represents the $i$-th component of vector $B_{\mu}(y)$, and $t_{1/2,\lambda\mu/(|y_{i}|+\epsilon_{i})^{1/2-p}}$ and $f_{1/2,\lambda\mu/(|y_{i}|+\epsilon_{i})^{1/2-p}}$ are obtained by replacing $\lambda$ with $\lambda\mu/(|y_{i}|+\epsilon_{i})^{1/2-p}$ in $t_{1/2,\lambda}$ and $f_{1/2,\lambda}$, respectively.
\end{lemma}
\textbf{proof.} We notice that, $\mathcal{C}^{1/2}_{\lambda,\mu}(z,y)$ can be rewritten as
\begin{eqnarray*}
\mathcal{C}^{1/2}_{\lambda,\mu}(z,y)&=&\|z-(y-\mu A^{T}Ay+\mu A^{T}b)\|_{2}^{2}+\lambda\mu \sum_{i=1}^{n}\frac{|z_{i}|^{1/2}}{(|y_{i}|+\epsilon_{i})^{1/2-p}}+\mu\|b\|_{2}^{2}\\
&&+\|y\|_{2}^{2}-\mu\|Ay\|_{2}^{2}-\|y-\mu A^{T}Ay+\mu A^{T}b\|_{2}^{2}\\
&=&\|z-B_{\mu}(y)\|_{2}^{2}+\lambda\mu \sum_{i=1}^{n}\frac{|z_{i}|^{1/2}}{(|y_{i}|+\epsilon_{i})^{1/2-p}}+\mu\|b\|_{2}^{2}+\|y\|_{2}^{2}-\mu\|Ay\|_{2}^{2}\\
&&-\|B_{\mu}(y)\|_{2}^{2}.
\end{eqnarray*}
This implies that minimizing $\mathcal{C}^{1/2}_{\lambda,\mu}(z,y)$ for any fixed $\lambda,\mu\in(0,\infty)$ and $y\in \mathbb{R}^{n}$ is equivalent to
\begin{equation}\label{equ23}
\min_{z\in \mathbb{R}^{n}}\Big\{\|z-B_{\mu}(y)\|_{2}^{2}+\lambda\mu \sum_{i=1}^{n}\frac{|z_{i}|^{1/2}}{(|y_{i}|+\epsilon_{i})^{1/2-p}}\Big\},
\end{equation}
i.e.,
\begin{equation}\label{equ24}
\min_{z\in \mathbb{R}^{n}}\Big\{\sum_{i=1}^{n}\Big((z_{i}-[B_{\mu}(y)]_{i})^{2}+\lambda\mu \frac{|z_{i}|^{1/2}}{(|y_{i}|+\epsilon_{i})^{1/2-p}}\Big)\Big\}.
\end{equation}
Noting that the summation of equation (\ref{equ25}) is separable; hence, solving equation (\ref{equ25}) is equivalent to solving the following $n$ subproblem, for
$i\in[1,2,\cdots,n]$,
\begin{equation}\label{equ25}
\min_{z\in \mathbb{R}^{n}}\Big\{(z_{i}-[B_{\mu}(y)]_{i})^{2}+\lambda\mu \frac{|z_{i}|^{1/2}}{(|y_{i}|+\epsilon_{i})^{1/2-p}}\Big\}.
\end{equation}
Therefore, the proof is completed by Lemma \ref{lem1}. $\hfill{} \Box$

\begin{theorem}\label{the1}
For any $\lambda,\mu\in(0,\infty)$, if $z^{*}=(z_{1}^{*}, z_{2}^{*},\cdots, z_{n}^{*})^{\top}$ is an optimal solution to the problem $(P_{p,\theta,\epsilon}^{\lambda})$ and $\mu$
satisfies $0<\mu<\frac{1}{\|A\|_{2}^{2}}$, then
\begin{equation}\label{equ26}
z_{i}^{\ast}=\left\{
    \begin{array}{ll}
      f_{1/2,\lambda\mu/(|z_{i}^{\ast}|+\epsilon_{i})^{1/2-p}}([B_{\mu}(z^{\ast})]_{i}), & \ \mathrm{if} \ {|[B_{\mu}(z^{\ast})]_{i}|>
      t_{1/2,\lambda\mu/(|z^{\ast}_{i}|+\epsilon_{i})^{1/2-p}};} \\
      0, & \ \mathrm{if} \ {|[B_{\mu}(z^{\ast})]_{i}|\leq t_{1/2,\lambda\mu/(|z^{\ast}_{i}|+\epsilon_{i})^{1/2-p}}.}
    \end{array}
  \right.
\end{equation}
\end{theorem}
\textbf{proof.} By condition $0<\mu<\frac{1}{\|A\|_{2}^{2}}$, we can get that
\begin{eqnarray*}
\mathcal{C}^{1/2}_{\lambda,\mu}(z,z^{*})&=&\mu\Big\{\|Az-b\|_{2}^{2}+\lambda \sum_{i=1}^{n}\frac{|z_{i}|^{1/2}}{(|z^{\ast}_{i}|+\epsilon_{i})^{1/2-p}}\Big\}\\
&&+\Big\{\|z-z^{*}\|_{2}^{2}-\mu\|Az-Az^{*}\|_{2}^{2}\Big\}\\
&\geq& \mu\Big\{\|Az-b\|_{2}^{2}+\lambda \sum_{i=1}^{n}\frac{|z_{i}|^{1/2}}{(|z^{\ast}_{i}|+\epsilon_{i})^{1/2-p}}\Big\} \\
&\geq&\mu\mathcal{C}^{1/2}_{\lambda}(z^{\ast})\\
&=& \mathcal{C}^{1/2}_{\lambda,\mu}(z^{*},z^{*})
\end{eqnarray*}
for any $z\in \mathbb{R}^{n}$. This implies that $z^{*}$ is a local minimizer of $\mathcal{C}^{1/2}_{\lambda,\mu}(z,z^{*})$ as long as $z^{*}$ is an optimal solution to the problem $(P_{p,1/2,\epsilon}^{\lambda})$. Combined with Lemma \ref{lem3}, we finish the proof.  $\hfill{} \Box$

Next, we present an iterative thresholding algorithm for solving the problem $(P_{p,1/2,\epsilon}^{\lambda})$ for all $p\in(0,1)$ based on the above theoretical analysis.

With the thresholding representation (\ref{equ26}), the $1/2-\epsilon$ algorithm for solving the regularization problem $(P_{p,1/2,\epsilon}^{\lambda})$ can be naturally defined as
\begin{equation}\label{equ27}
z_{i}^{k+1}=h_{1/2,\lambda\mu/(|z_{i}^{k}|+\epsilon_{i})^{1/2-p}}([B_{\mu}(z^{k})]_{i}),\ \ \ k=0,1,2,\cdots,
\end{equation}
where $B_{\mu}(z^{k})=z^{k}+\mu A^{\top}(b-Az^{k})$, and $h_{1/2,\lambda\mu/(|z_{i}^{k}|+\epsilon_{i})^{1/2-p}}$ is obtained by replacing $\lambda$ with $\lambda\mu/(|z_{i}^{k}|+\epsilon_{i})^{1/2-p}$ in $h_{1/2,\lambda}$.

In general, the quality of the solution to a regularization problem depends seriously on the setting of the regularization parameter $\lambda>0$.
Suppose that the vector $z^{\ast}$ of sparsity $r$ is the optimal solution of the regularization problem $(P_{p,1/2,\epsilon}^{\lambda})$.
In $1/2-\epsilon$ algorithm, we set
\begin{equation}\label{equ37}
\lambda=\frac{8(\lceil B_{\mu}(z^{k})\rfloor_{r+1})^{3/2}(\lceil z^{k}\rfloor_{r+1}+\lceil\epsilon\rfloor_{r+1})^{1/2-p}}{\sqrt{54}\mu}
\end{equation}
in each iteration, where $\lceil\rfloor$ is defined in Definition \ref{def1}, and $\lceil B_{\mu}(z^{k})\rfloor_{i}$ ($\lceil z^{k}\rfloor_{i}$, $\lceil\epsilon\rfloor_{i}$) represents the $i$-th component of vector $\lceil B_{\mu}(z^{k})\rfloor$ ($\lceil z^{k}\rfloor$, $\lceil\epsilon\rfloor$) for all $i\in[1,2,\cdots,n]$. When doing so, the $1/2-\epsilon$ algorithm will be adaptive and free from the choice of regularization parameter.

\begin{algorithm}
\caption{: The $1/2-\epsilon$ algorithm}
\label{alg:A}
\begin{algorithmic}
\STATE {\textbf{Initialize}: Choose $x^{0}\in \mathbb{R}^{n}$, $\epsilon_{i}>0$, $\lambda_{0}>0$, $\mu_{0}=\frac{1-\eta}{\|A\|_{2}^{2}}(\eta\in(0,1))$ and $p\in(0,1)$;}
\STATE {$k=0$;}
\STATE {\textbf{while} not converged \textbf{do}}
\STATE \ \ \ \ \ \ \ {$B_{\mu}(z^{k})=z^{k}+\mu A^{\top}(b-Az^{k})$;}
\STATE \ \ \ \ \ \ \ {$\lambda=\frac{8(\lceil B_{\mu}(z^{k})\rfloor_{r+1})^{3/2}(\lceil z^{k}\rfloor_{r+1}+\lceil\epsilon\rfloor_{r+1})^{1/2-p}}{\sqrt{54}\mu}$, $\mu=\mu_{0}$;}
\STATE \ \ \ \ \ \ \ \ {$t=t_{1/2,\lambda\mu/(|z_{i}^{k}|+\epsilon_{i})^{1/2-p}}=\frac{\sqrt[3]{54}}{4}(\lambda\mu/(|z_{i}^{k}|+\epsilon_{i})^{1/2-p})^{2/3}$;}
\STATE \ \ \ \ \ \ \ \ {for\ $i=1:n$}
\STATE \ \ \ \ \ \ \ \ \ \ \ \ \ {1.\ $|[B_{\mu}(z^{k}]_{i}|>t$, then $z^{k+1}_{i}=f_{1/2,\lambda\mu/(|z_{i}^{k}|+\epsilon_{i})^{1/2-p}}([B_{\mu}(z^{k}]_{i})$}
\STATE \ \ \ \ \ \ \ \ \ \ \ \ \ {2.\ $|[B_{\mu}(z^{k}]_{i}|\leq t$, then $z^{k+1}_{i}=0$}
\STATE {$k\rightarrow k+1$;}
\STATE{\textbf{end while}}
\STATE{\textbf{return}: $z^{k+1}$}
\end{algorithmic}
\end{algorithm}

\begin{remark}\label{re1}
It is worth emphasizing that the $1/2-\epsilon$ algorithm reduces to the Half algorithm\cite{xu2} when we set $p=1/2$.
\end{remark}

In the following, we provide some convergence analysis for the $1/2-\epsilon$ algorithm under some specific conditions.

\begin{theorem}\label{the2}
Let $\{z^{k}\}$ be the sequence generated by the $1/2-\epsilon$ algorithm with the step size $\mu$ satisfying $0<\mu<\frac{1}{\|A\|_{2}^{2}}$. Then
\begin{description}
\item[$\mathrm{1)}$] The sequence $\{z^{k}\}$ is a minimization sequence, and the sequence $\{\mathcal{C}^{1/2}_{\lambda}(z^{k})\}$ converging to
$\mathcal{C}^{1/2}_{\lambda}(z^{\ast})$, where $z^{\ast}$ is a limit point of minimization sequence $\{z^{k}\}$;
\item[$\mathrm{2)}$] The sequence $\{z^{k}\}$ is asymptotically regular, i.e., $\lim_{k\rightarrow\infty}\|z^{k+1}-z^{k}\|_{2}^{2}=0$;
\item[$\mathrm{3)}$] Any accumulation point of the sequence $\{z^{k}\}$ is a stationary point of the problem $(P_{p,1/2,\epsilon}^{\lambda})$.
\end{description}
\end{theorem}
\textbf{proof.} The proof of above theorem follows from the fact that the step size $\mu$ satisfying $0<\mu<\frac{1}{\|A\|_{2}^{2}}$ and a similar argument as used in the
proof of [2, Theorem 3]. $\hfill{} \Box$

\subsection{The $2/3-\epsilon$ algorithm for solving the problem $(P_{p,2/3,\epsilon}^{\lambda})$}\label{subsection3-2}
In the subsection, we propose the $2/3-\epsilon$ algorithm to solve the problem $(P_{p,2/3,\epsilon}^{\lambda})$ for all $p\in(0,1)$.

For any $\lambda,\mu\in(0,\infty)$, $p\in(0,1)$ and $z, y\in \mathbb{R}^{n}$, let
\begin{equation}\label{equ29}
\mathcal{C}^{2/3}_{\lambda}(z)=\|Az-b\|_{2}^{2}+\lambda \sum_{i=1}^{n}\frac{|z_{i}|^{2/3}}{(|z_{i}|+\epsilon_{i})^{2/3-p}}
\end{equation}
and
\begin{equation}\label{equ30}
\begin{array}{llll}
\mathcal{C}_{\lambda,\mu}^{2/3}(z,y)&=&\mu\|Az-b\|_{2}^{2}+\displaystyle\lambda\mu \sum_{i=1}^{n}\frac{|z_{i}|^{2/3}}{(|y_{i}|+\epsilon_{i})^{2/3-p}}-\mu\|Az-Ay\|_{2}^{2}+\|z-y\|_{2}^{2}.
\end{array}
\end{equation}
Similar argument as the generation of $1/2-\epsilon$ algorithm, the $2/3-\epsilon$ algorithm for solving the problem $(P_{p,2/3,\epsilon}^{\lambda})$ can be defined as

\begin{equation}\label{equ31}
z_{i}^{k+1}=h_{2/3,\lambda\mu/(|z_{i}^{k}|+\epsilon_{i})^{2/3-p}}([B_{\mu}(z^{k})]_{i}),\ \ \ k=0,1,2,\cdots,
\end{equation}
where $B_{\mu}(z^{k})=z^{k}+\mu A^{\top}(b-Az^{k})$, and $h_{2/3,\lambda\mu/(|z_{i}^{k}|+\epsilon_{i})^{2/3-p}}$ is obtained by replacing $\lambda$ with $\lambda\mu/(|z_{i}^{k}|+\epsilon_{i})^{2/3-p}$ in $h_{2/3,\lambda}$.

\begin{algorithm}
\caption{: The $2/3-\epsilon$ algorithm}
\label{alg:B}
\begin{algorithmic}
\STATE {\textbf{Initialize}: Choose $x^{0}\in \mathbb{R}^{n}$, $\epsilon_{i}>0$, $\mu_{0}=\frac{1-\eta}{\|A\|_{2}^{2}}(\eta\in(0,1))$ and $p\in(0,1)$;}
\STATE {$k=0$;}
\STATE {\textbf{while} not converged \textbf{do}}
\STATE \ \ \ \ \ \ \ {$B_{\mu}(z^{k})=z^{k}+\mu A^{\top}(b-Az^{k})$;}
\STATE \ \ \ \ \ \ \ {$\lambda=\frac{\sqrt[3]{4^{4}}(\lceil B_{\mu}(z^{k})\rfloor_{r+1})^{4/3}(\lceil z^{k}\rfloor_{r+1}+\lceil\epsilon\rfloor_{r+1})^{2/3-p}}{\sqrt[9]{48^{4}}\mu}$, $\mu=\mu_{0}$;}
\STATE \ \ \ \ \ \ \ \ {$t=t_{2/3,\lambda\mu/(|z_{i}^{k}|+\epsilon_{i})^{2/3-p}}=\frac{\sqrt[4]{48}}{3}(\lambda\mu/(|z_{i}^{k}|+\epsilon_{i})^{2/3-p})^{3/4}$;}
\STATE \ \ \ \ \ \ \ \ {for\ $i=1:n$}
\STATE \ \ \ \ \ \ \ \ \ \ \ \ \ {1.\ $|[B_{\mu}(z^{k}]_{i}|>t$, then $z^{k+1}_{i}=f_{2/3,\lambda\mu/(|z_{i}^{k}|+\epsilon_{i})^{2/3-p}}([B_{\mu}(z^{k}]_{i})$}
\STATE \ \ \ \ \ \ \ \ \ \ \ \ \ {2.\ $|[B_{\mu}(z^{k}]_{i}|\leq t$, then $z^{k+1}_{i}=0$}
\STATE {$k\rightarrow k+1$;}
\STATE{\textbf{end while}}
\STATE{\textbf{return}: $z^{k+1}$}
\end{algorithmic}
\end{algorithm}

In $2/3-\epsilon$ algorithm, we set the regularization parameter $\lambda$ as
\begin{equation}\label{equ32}
\lambda=\frac{\sqrt[3]{4^{4}}(\lceil B_{\mu}(z^{k})\rfloor_{r+1})^{4/3}(\lceil z^{k}\rfloor_{r+1}+\lceil\epsilon\rfloor_{r+1})^{2/3-p}}{\sqrt[9]{48^{4}}\mu}
\end{equation}
in each iteration.

Similar argument as the Theorem \ref{the2}, the sequence $\{X^{k}\}$ generated by the $2/3-\epsilon$ algorithm is a minimization sequence and asymptotically regular.
Moreover, any accumulation point of $\{X^{k}\}$ is a stationary point of the problem $(P_{p,2/3,\epsilon}^{\lambda})$. Its proof also follows from the fact that the
step size $\mu$ satisfying $0<\mu<\frac{1}{\|A\|_{2}^{2}}$ and a similar argument as used in the proof of [2, Theorem 3].

\begin{remark}\label{re2}
If we set $p=2/3$, the $2/3-\epsilon$ algorithm reduces to the $2/3$ algorithm \cite{cao3}.
\end{remark}

\section{Numerical experiments}\label{section4}
In this section, we first present numerical results of the $1/2-\epsilon$ algorithm and $2/3-\epsilon$ algorithm for some sparse signal recovery problems and then
compare them with some state-of-art methods including iterative hard thresholding algorithm (Hard algorithm)\cite{thom6}, iterative soft thresholding algorithm (Soft
algorithm)\cite{dau7}, Half algorithm \cite{xu2} and 2/3 algorithm \cite{cao3} in some sparse recovery problems. We generate a measurement matrix $A\in \mathbb{R}^{m\times n}$
with entries independently drawn by random from a Gaussian distribution, $\mathcal{N}(0,1)$. To show the success rate of these algorithms in the recovery of the sparse
signals with the different sparsity for a given measurement matrix $A\in \mathbb{R}^{m\times n}$, we randomly generate sparse vectors $z_{0}\in \mathbb{R}^{n}$ and generate
vectors $b$ by $b=Az_{0}$. Therefore, we know the sparsest solution to the linear system $b=Az_{0}$. In our experiments, the stopping criterion is defined as
$$\frac{\|z^{k+1}-z^{k}\|_{2}}{\|z^{k}\|_{2}}\leq \mathrm{Tol}$$
where $z^{k+1}$ and $z^{k}$ are numerical results from two continuous iterative steps and $\mathrm{Tol}$ is a small given number (we set $\mathrm{Tol}=10^{-8}$).
The success is measured by computing the relative error (RE):
$$\mathrm{RE}=\frac{\|z^{\ast}-z_{0}\|_{2}}{\|z_{0}\|_{2}}$$
to indicate a perfect recovery of the original sparse vector $z_{0}$, and the success is declared when $\mathrm{RE}\leq 10^{-4}$. Moreover, we adapt $\epsilon$ at each
iteration as a function of the current guess $z^{k}$, and set
$$\epsilon_{i}=\max\{\gamma|[\mu A^{\top}(b-Az^{k})]_{i}|, 10^{-3}\}, \ \ \ \gamma=0.7.$$
In all of our experiments, we repeat 20 tests and present average results.
The experiments are all performed on a Lenovo-PC with an Intel(R) Core(TM) i7-6700 CPU @3.40GHZ
with 16GB of RAM running Microsoft Windows 7.

\subsection{Performance of $1/2-\epsilon$ algorithm and $2/3-\epsilon$ algorithm}\label{subsection4-1}
In this subsection, we carry out a series of experiments to demonstrate the performance of the $1/2-\epsilon$ algorithm and the $2/3-\epsilon$ algorithm for some sparse
signal recovery problems. In our experiments, we set $m=128$, $n=512$.

\begin{figure}[h!]
 \centering
 \includegraphics[width=0.6\textwidth]{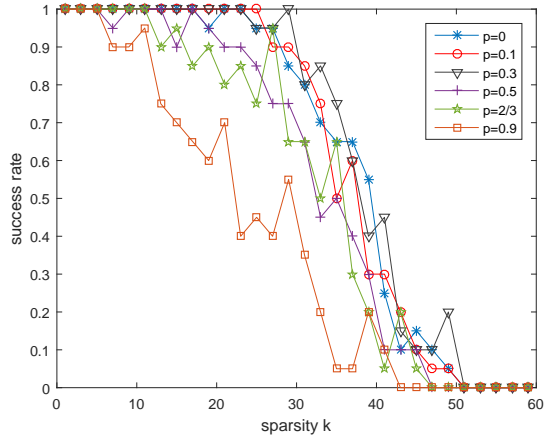}
\caption{The comparison of the success rate for the $1/2-\epsilon$ algorithm in the recovery of a sparse signal with different parameter values $p$.}
\label{fig1}       
\end{figure}

\begin{figure}[h!]
 \centering
 \includegraphics[width=0.6\textwidth]{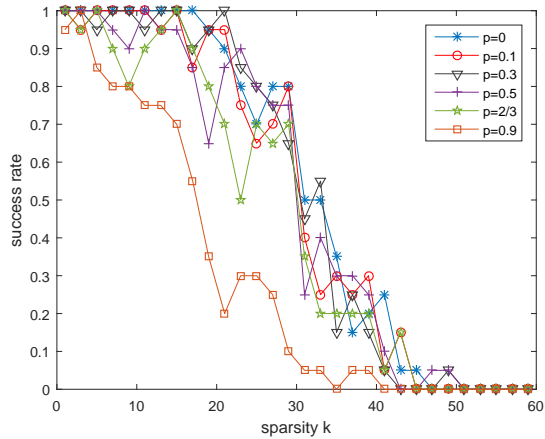}
\caption{The comparison of the success rate for the $2/3-\epsilon$ algorithm in the recovery of a sparse signal with different parameter values $p$.}
\label{fig2}       
\end{figure}

The graphs presented in Figure \ref{fig1} and Figure \ref{fig2} show the success rate of $1/2-\epsilon$ algorithm and $2/3-\epsilon$ algorithm in recovering the true (sparsest)
solution. From Figures \ref{fig1} and \ref{fig2}, we can see that $1/2-\epsilon$ algorithm can exactly recover the ideal signal until $k$ is around $25$ when $p=0.1$, and $2/3-\epsilon$
algorithm's counterpart is around $17$ when $p=0$. As we can see, the parameter $p=0.1$ is the best strategy for the $1/2-\epsilon$ algorithm, and the the parameter $p=0$ is the best strategy
for the $2/3-\epsilon$ algorithm.

\subsection{Compared with some state-of-art methods}\label{subsection4-2}
\begin{figure}[h!]
 \centering
 \includegraphics[width=0.6\textwidth]{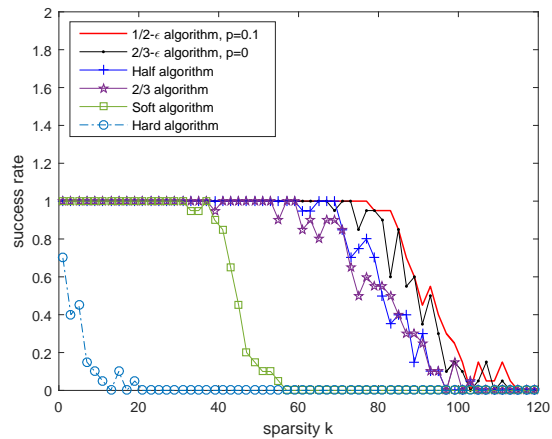}
\caption{The success rate of six algorithms in the recovery of a sparse signal with different sparsity (noiseless case).}
\label{fig3}       
\end{figure}

\begin{figure}[h!]
 \centering
 \includegraphics[width=0.6\textwidth]{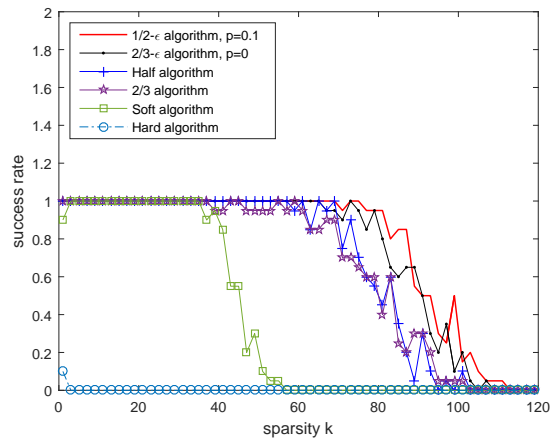}
\caption{The success rate of six algorithms in the recovery of a sparse signal with different sparsity (noise case).}
\label{fig4}       
\end{figure}

In this subsection, we compare our algorithms (the $1/2-\epsilon$ algorithm and $2/3-\epsilon$ algorithm) with some state-of-art methods including Hard
algorithm\cite{thom6}, Soft algorithm\cite{dau7}, Half algorithm\cite{xu2} and $2/3$ algorithm\cite{cao3} in some sparse recovery problems.
In our experiments, we set $p=0.001$ in the $1/2-\epsilon$ algorithm and the $2/3-\epsilon$ algorithm, and set $m=256$, $n=1024$ to size the dimension
of the matrix $A\in \mathbb{R}^{m\times n}$ and the length of the vector $z_{0}\in \mathbb{R}^{n}$. Two different cases in sparse recovery problems will
be considered: exactly sparse signals recovery in the noiseless case and exactly sparse signals recovery in the noise case. In noiseless case, we generate
vectors $b$ by $b=Az_{0}$, where $A\in \mathbb{R}^{256\times 1024}$ is a measurement matrix with entries independently drawn by random from a Gaussian
distribution, $\mathcal{N}(0,1)$, and $z_{0}\in \mathbb{R}^{1024}$ is a randomly sparse vector $z_{0}\in \mathbb{R}^{1024}$. Turn to the noise case,
we use the same matrix $A$, and generate a random vector $z_{0}$ with a prespecified cardinality of nonzeros. We compute $b=Az_{0}+e$, where
$e\in \mathcal{N}(0,\sigma)$ $(\sigma=10^{-5})$. Thus, the original vector $z_{0}$ is a feasible solution and close to the optimal solution. Due to
the presence of noise, it becomes harder to accurately recover the original signal $z_{0}$.
The graphs presented in Figure \ref{fig3} and Figure \ref{fig4} show the success rate of $1/2-\epsilon$ algorithm ($p=0.1$), $2/3-\epsilon$ algorithm ($p=0$),
Half algorithm, $2/3$ algorithm, Soft algorithm and Hard algorithm in recovering the true (sparsest) solution. From Figure \ref{fig3}, we can see that the $1/2-\epsilon$
algorithm ($p=0.1$) can exactly recover the ideal signal until $k$ is around $78$, and the $2/3-\epsilon$ algorithm ($p=0$) is around $70$. The results in
noise case are consistent with the noiseless case. We can see that the $1/2-\epsilon$ algorithm ($p=0.1$) again has the best performance in recovering the sparse signals
in the six algorithms with noise or not, and the $2/3-\epsilon$ algorithm ($p=0$) performs the second best.

\section{Conclusions}\label{section5}

In this paper, we studied two modified $\l_{p}$-norm regularization minimization problems $(P_{p,1/2,\epsilon}^{\lambda})$ and $(P_{p,2/3,\epsilon}^{\lambda})$ to approximate
the NP-hard problem $(P_{p}^{\lambda})$ for all $p\in(0,1)$. Inspired by the good performances of iterative $\l_{1/2}$ thresholding algorithm and iterative $\l_{2/3}$ algorithm
in some sparse signal recovery problems, the $1/2-\epsilon$ algorithm and $2/3-\epsilon$ algorithm are generated to solve the problems $(P_{p,1/2,\epsilon}^{\lambda})$  and 
$(P_{p,2/3,\epsilon}^{\lambda})$ for all $p\in(0,1)$. Numerical results show that our algorithms perform effectively in finding the sparse signals in some sparse signal recovery 
problems for some proper $p\in(0,1)$. Moreover, the numerical results also show that the $1/2-\epsilon$ algorithm performs the best in some sparse signal recovery problems compared 
with some state-of-art methods, and the $2/3-\epsilon$ algorithm performs the second best for some proper $p$.

\section*{Acknowledgments}
We would like to thank editors and reviewers for their comments which help us to enrich the content and improve the presentation of the results in this paper.
The work was supported by the National Natural Science Foundations of China (11771347, 91730306, 41390454, 11271297) and the Science Foundations of Shaanxi Province
of China (2016JQ1029, 2015JM1012).




\section*{References}

\begin{thebibliography}{00}



\bibitem{chen1}
X. Chen, D. Ge, Z. Wang, and Y. Ye, Complexity of unconstrained $L_{2}-L_{p}$ minimization. Mathematical Programming, 143 (2014) 371--383.


\bibitem{xu2}
Z. Xu, X. Chang, F. Xu, and H, Zhang. L1/2 Regularization: A thresholding representation theory and a fast solver. IEEE Transactions on Neural Networks and Learning Systems, 24(7) (2012) 1013--1027.

\bibitem{cao3}
W. Cao, J. Sun, Z. Xu, Fast image deconvolution using closed-form thresholding formulas of $L_{q}(q=\frac{1}{2},\frac{2}{3})$ regularization. Journal of Visual Communication and Image Representation,
24(1) (2013) 31--41.

\bibitem{chen4}
X. Chen, F. Xu, Y. Ye, Lower bound theory of nonzero entries in solutions of $\ell_{2}-\ell_{p}$ minimization. SIAM Journal on Scientific Computing, 32(5) (2010) 2832--2852.

\bibitem{chen5}
X. Chen, L. Niu, Y. Yuan, Optimality conditions and a smoothing trust region newton method for nonlipschitz optimization. SIAM Journal on Optimization, 23(3) (2013) 1528--1552.

\bibitem{thom6}
T. Blumensath, M. E. Davies, Iterative thresholding for sparse approximations. Journal of Fourier Analysis and Applications, 14(5-6) (2008) 629--654.

\bibitem{dau7}
I. Daubechies, M. Defrise, C. De Mol, An iterative thresholding algorithm for linear inverse problems with a sparsity constraint. Communications on Pure and Applied Mathematics, 57 (11) (2004)
1413--1457.

\bibitem{don8}
D. L. Donoho, Denoising by soft-thresholding. IEEE Transactions on Information Theory, 41(3) (1995) 613--627.

\bibitem{tom9}
T. Goldstein, S. Osher, The split Bregman method for L1-regularized problems. SIAM Journal on Imaging Sciences, 2(2) (2009) 323--343.

\bibitem{yin10}
W. Yin, S. Osher, D. Goldfarb, and J. Darbon, Bregman iterative algorithms for $\ell_{1}$-minimization with applications to compressed sensing. SIAM Journal on Imaging Sciences, 1(1) (2008) 143--168.

\bibitem{yang11}
J. Yang, Y. Zhang, Alternating direction algorithms for $\ell_{1}$ problems in compressive sensing. SIAM Journal on Scientific Computing, 33(1) (2011) 250--278.

\bibitem{sim12}
S. Foucart, H. Rauhut, A mathematical introduction to compressive sensing. Springer, New York (2010).



\end{thebibliography}


\end{document}